\documentclass{amsart}
\usepackage{graphics,graphicx} 

\usepackage[T1]{fontenc}

\usepackage[cp1250]{inputenc}

\usepackage[arrow, matrix, curve]{xy}
\usepackage{amsmath,amsthm,amssymb}
\usepackage{amscd}

\newenvironment{Proof of}[1]{\emph{Proof of #1.}}{$\qquad \square$\par}

\DeclareMathOperator{\clsp}{\overline{span}}

\renewcommand{\a}{\widetilde a}

\renewcommand{\b}{\widetilde b}

\newcommand{\X}{\widehat{X}}

\newcommand{\tal}{\widetilde \alpha}

\newcommand{\al}{\alpha}

\newcommand{\A}{\mathcal A}

\newcommand{\B}{\mathcal B}

\newcommand{\bbZ}{\mathbb Z}
\newcommand{\bbN}{\mathbb N}

\newtheorem{theorem}{Theorem}[section] 
\newtheorem{proposition}[theorem]{Proposition} 

\theoremstyle{definition} 
\newtheorem{definition}[theorem]{Definition}

\newtheorem{remark}[theorem]{Remark}



\author{Bartosz Kosma  Kwa\'sniewski}
\address{ Institute of Mathematics, Polish Academy of Science,  ul. \'Sniadeckich 8, PL-00-956 Warszawa, Poland // Institute of Mathematics, University  of Bialystok,  ul. Akademicka 2, PL-15-267  Bialystok, Poland}
\email{bartoszk@math.uwb.edu.pl}

\title{Extensions of $C^*$-dynamical systems  to systems  with \\ complete transfer operators}

\begin{document}
\maketitle

\begin{abstract}
Starting from an arbitrary endomorphism $\alpha$ of a
unital $C^*$-algebra $\A$ we construct a bigger algebra $\B$
and extend $\alpha$ onto $\B$ in such a way that $\alpha:\B \to \B$
has a unital kernel and a hereditary range, so there exists a
unique non-degenerate transfer operator for $(\B,\alpha)$, called the complete transfer operator. The pair  $(\B,\alpha)$ is universal with respect to a suitable notion of a covariant representation and depends on a choice of an ideal in $\A$.  The construction enables a natural definition of the crossed product  for arbitrary $\alpha$.
\\
\end{abstract}


\section*{Introduction}


The crossed-product of a unital $C^*$-algebra $\A$ by an automorphism
$\alpha:\A\to \A$ is  the universal $C^*$-algebra
generated by a copy of $\A$ and a unitary element $U$ satisfying
the relations
$$
\alpha(a)= Ua U^*,\qquad \alpha^{-1}(a)=U^*aU, \qquad a\in \A.
$$
Algebras arising in this way (or their versions adapted to actions of groups of auto\-mor\-phisms) are very well understood and are one of the standard constructions in $C^*$-theory.
 On the other hand, the natural desire to adapt this kind of constructions to endomorphisms (or semigroups of endomorphisms)  encounters serious obstacles from the very beginning. Roughly, it is  caused by   the irreversibility  of the system $(\A,\alpha)$   - the lack of the map $\alpha^{-1}$.

The difficulty of the matter manifests itself  in a variety of
approaches, see, for example,  \cite{CK}, \cite{Paschke},
\cite{Stacey}, \cite{Murphy},  \cite{exel1}, \cite{exel2},
 \cite{Ant-Bakht-Leb},  \cite{kwa}, which do however have a certain
nontrivial intersection. They  mostly agree, and simultaneously
boast their greatest successes, in the  case when the dynamics is
implemented by a monomorphism with a hereditary range. In view
of  \cite{Bakht-Leb},
\cite{Ant-Bakht-Leb},  see also \cite{kwa3}, this coincidence seems to be well
understood. Namely, it was noticed in \cite{Bakht-Leb},
\cite{Ant-Bakht-Leb} that for a class of endomorphisms (as shown in \cite{kwa3} consisting of endomorphisms with a complementary kernel and a hereditary range)  there exists a unique non-degenerate transfer operator  \cite{exel2},
called by the authors of  \cite{Bakht-Leb} the \emph{complete transfer
operator}. In this case the theory goes smooth, in the  spirit very similar to that of crossed products by automorphisms,    as the complete transfer operator  $\mathcal{L}$ takes
over the role classically played by $\alpha^{-1}$. The authors of  \cite{Ant-Bakht-Leb}  showed that all of the construction in the aforementioned various approaches can be reduced to the crossed product for systems $(\A,\alpha)$ with  complete transfer operators. At the end of \cite{Ant-Bakht-Leb} they  argued that a general crossed product construction should consist of two steps:
\begin{itemize}
\item[1)] ``an initial object and an extension procedure'';
\item[2)] ``the crossed product for systems with complete transfer operators''.
\end{itemize}
Our   goal  is to provide the missing first step in the  above   scheme  when the initial object is an arbitrary endomorphism of a unital $C^*$-algebra.
The previous preprint version of the present note together with \cite{Ant-Bakht-Leb}  enabled the authors of  \cite{kwa-leb1}  to develop a general approach to crossed products by arbitrary endomorphisms. We also hope to use the elaborated  $C^*$-algebraic extension method in the further more detailed  analysis of $C^*$-algebras of this type, cf., e.g., \cite[problem on p. 1830]{kwa-leb1}. Such an approach already proved to be very useful in the case the initial system is commutative, see  \cite{kwa}, \cite{maxid}, \cite{kwa2} and  Remark \ref{commutative remark} below.  

We note that the idea behind  our construction is very similar to that of known as dilation, see for instance  \cite{La} and references therein. However, the dilation of endomor\-phisms to automorphisms, one can come across the in literature, applies to injective endo\-morphisms and yield crossed products that are equivalent only up to the Morita equi\-valence. Whilst our extension    procedure concern  general (not necessarily injective) endomorphisms and the crossed product for  the initial and the extended object are naturally  isomorphic, see Theorem \ref{rozszerzenie a repr. kowariantne2} ii) below.

\section{The problem and the spatial operator-algebraic considerations}\label{reprezentacje kowariantne}
Throughout this paper we let  $\A$ be a $C^*$-algebra with an identity $1$. By a \emph{$C^*$-dynamical system} we mean a pair $(\A,\alpha)$ where  $\alpha:\A\to \A$ is an endomorphism of $\A$ (all the morphisms appearing in the text are assumed to be $^*$-preserving).
To explain the problem let us suppose that the system $(\A,\alpha)$ is faithfully represented on a Hilbert space $
H$. We assume  $\A$ is a $C^*$-subalgebra of the algebra $L(H)$ of all bounded linear operators on $H$,  $1$ is the identity operator, and there is  $U\in L(H)$ such that 
$$
\alpha(a)=UaU^*, \qquad a\in \A.
$$ 
By   \cite[Prop. 2.2]{Leb-Odz} and \cite[Lem 1.2]{kwa-leb1} the multiplicativity of $\alpha$ is equivalent to requiring that 
$$
U\textrm{ is a (power) partial isometry and } U^*U \in\A',
$$
where  $\A'$ denotes the commutant of $\A$. Hence by \cite[Prop. 2.2. and Prop 3.10]{Leb-Odz} 
$$
\B:=\clsp\{U^{*n} a U^{n}: a\in \A, \,\, n \in \mathbb{N}\},
$$
is a  minimal $C^*$-algebra containing $\A$ and such that
  the following relations hold
 $$
 U\B U^*\subset \B, \qquad U^*\B U\subset \B,\qquad U^*U\in Z(\B)=\B'\cap \B.
  $$
  Therefore, see also  \cite[3.1]{Bakht-Leb} or \cite[2.5]{Ant-Bakht-Leb}, putting 
  $$
  \alpha(b):=Ub U^*, \qquad \mathcal{L}(b):=U^* b U, \qquad b\in \B,
  $$
  we obtain an endomorphism  $\alpha:\B\to \B$  that extends $\alpha:\A\to \A$,\footnote{as a rule we use the same symbol for  endomorphisms and their extensions} and a linear  operator $\mathcal{L}:\B\to \B$ which is a complete transfer operator for the extended system $(\B,\alpha)$  (we recall the definition of the complete transfer operator in section \ref{transfer}). In the present note we give a positive answer to the following 
\medskip

\begin{quote}
\textbf{Question:}	Does there exist an efficient  description  of the triple $(\B,\alpha, \mathcal{L})$ in terms of the initial $C^*$-dynamical system $(\A,\alpha)$, independent of the representation in $L(H)$? 
\end{quote}
 \begin{remark}\label{commutative remark}
By  \cite[Prop. 4.1]{Leb-Odz}, if $\A$ is commutative  then $\B$ is also commutative. Hence in this case the $C^*$-dynamical systems
 $(\A,\alpha)$  and $(\B,\alpha)$ correspond to topological dynamical systems $(X,\al)$
 and $(\X,\tal)$ respectively (consisting of  compact Hausdorff spaces and partial mappings). The description of  $(\X,\tal)$ only in terms of $(X,\al)$ was obtained in \cite{maxid} under the additional assumption that $U^*U\in \A$. In general,  $(\X,\tal)$ is described by  $(X,\al)$  in \cite{kwa2} but this description requires additional data encoded in the ideal
 \begin{equation}\label{ideal J}
 J:=U^*U\A\cap \A=\{a\in \A: U^*Ua=a\},
\end{equation}
 cf. also \cite{kwa-leb1}.  The relationship between $(X,\al)$ and $(\X,\tal)$ is of particular interest.
In   \cite{kwa2}, \cite{maxid}, \cite{kwa}, a number of examples is studied and it is shown  that  the space $\X$ may be viewed as a generalization of the topological
  inverse limit space and as a rule $\X$ is topologically very complicated - contains
   \emph{indecomposable continua}, has a  structure of \emph{hyperbolic attractors},
    or of a space arising from \emph{substitution tilings}.  Thus, among the other things, the  construction of the present paper could be considered as a tool to obtain  non-commutative counterparts of the aforementioned objects.
    \end{remark}
We will analyze the structure of $\B$ by means of the following 'approximating' algebras 
$$
\label{algebry B_n 1}
\B_n:=\overline{\{\sum_{i=0}^n U^{*i}a_i U^i: a_i \in \A, i=0,...,n\}}, \qquad n\in \bbN,
$$
see  \cite[Prop. 3.8 (ii)]{Leb-Odz}. The family  $\{\B_n\}_{n\in\bbN}$ fixes  the structure of a direct limit on $\B$:
$$
\A=\B_0\subset \B_1 \subset ...\subset \B_n\subset ...,\qquad  \textrm{ and }\qquad \B=\overline{\bigcup_{n\in\bbN} \B_n}.
$$
The first crucial  step is to notice that the algebras $\B_n$ can  be canonically identified with direct sums of subalgebras of
the $C^*$-algebra $C^*(\A,U^*U)$ generated by $\A$ and $U^*U$. 
\begin{proposition}\label{kowariantna podpowiedz do algebry wpsolczynnikow}
Let $n\in \mathbb{N}$. Every element
$a\in \B_n$ can be presented in the form
$$
a = a_0 +U^*a_1U +... + U^{*n}a_n U^n
$$
where
\begin{equation}\label{coefficients warunki}
a_i\in (1-U^*U)\alpha^i(1)\A \alpha^i(1), \quad  i=0,...,n-1,\qquad a_n \in \alpha^n(1)\A \alpha^n(1),
\end{equation}
and this form is unique. Actually,
$a\mapsto a_0\oplus a_1\oplus ...\oplus a_n$
establishes the  isomorphism
\begin{equation}\label{coefficients izomor muchomor}
\B_n\cong (1-U^*U)\A \oplus (1-U^*U)\alpha^1(1)\A \alpha^1(1)\oplus ...\oplus \alpha^n(1)\A \alpha^n(1).
\end{equation}
\end{proposition}

\textbf{Proof.}
Let $a\in \B_n$. Then $a=\sum_{i=0}^{n}\mathcal{L}^i(b_i)$ where  $b_i\in \A$ and
$\mathcal{L}(\cdot)=U^*(\cdot)U$. Without loss of generality we may assume
that  $b_i\in \alpha^i(1)\A\alpha^i(1)$, because $\mathcal{L}^i( b_i )=\mathcal{L}^i(\alpha^i(1) b_i\alpha^i(1))$. We recall, cf. \cite[Prop. 3.6 (iv)]{Leb-Odz}, that 
the family $\{\mathcal{L}^k(1)\}_{k\in \bbN} \subset Z(\B)$ is
a decreasing  sequence of orthogonal projections.
 We will construct elements  $a_i$ satisfying \eqref{coefficients warunki} modifying inductively the elements $b_i$.
For  $a_0$ we take  $b_0(1- \mathcal{L}(1))$  and  'the remaining part'  of $b_0$
we include in $b_1$, that is  we put $c_1=b_1 + \alpha(b_0)$. Then  $a= a_0+ \mathcal{L}(c_1)+... +\mathcal{L}^n(b_n)$, because $b_0\mathcal{L}(1)=\mathcal{L}(\alpha(b_0))$.
\\
Continuing in this manner we get $k<n$ coefficients
$a_0,...,a_{k-1}$ satisfying \eqref{coefficients warunki} and such
that
$a=a_0+...+\mathcal{L}^{k-1}(a_{k-1})+\mathcal{L}^{k}(c_{k})+\mathcal{L}^{k+1}(b_{k+1}) +...
+\mathcal{L}^n(b_n)$ and $c_k\in\alpha^k(1)\A\alpha^k(1)$.  We put $a_k=c_k(1- \mathcal{L}(1))\in \A$ and
$c_{k+1}=b_{k+1} +\alpha(c_k)$. Then  $a_k
\in(1-U^*U)\alpha^k(1)\A\alpha^k(1)$ and  the following
computations
\begin{align*}
\mathcal{L}^k(c_k)&=\mathcal{L}^k(c_k)\mathcal{L}^k(1)=\mathcal{L}^k(c_k)\big(\mathcal{L}^k(1)-\mathcal{L}^{k+1}(1)\big)
+ \mathcal{L}^{k+1}(1)\mathcal{L}^k(c_k)
\\
&=\mathcal{L}^k\big(c_k (1-\mathcal{L}(1))\big)+ \mathcal{L}^{k+1}(\alpha(c_k))=\mathcal{L}^k(a_k)+ \mathcal{L}^{k+1}(\alpha(c_k))
\end{align*}
 show that
$a=a_0+...+\mathcal{L}^{k}(a_{k})+\mathcal{L}^{k+1}(c_{k+1})+...
+\mathcal{L}^n(b_n)$.
\\
Thus we may assume that  \eqref{coefficients warunki} holds. These conditions  imply that
$$
\mathcal{L}^i(a_i)\in \big(\mathcal{L}^i(1)-\mathcal{L}^{i+1}(1)\big)\A, \qquad i=0,...,n-1.
$$
Since $\{\mathcal{L}^k(1)\}_{k\in \bbN} \subset Z(\B)$ are decreasing  orthogonal projections, the projections  $1-\mathcal{L}(1)$, $\mathcal{L}(1)-\mathcal{L}^2( 1)$,\,...,
$\mathcal{L}^{n}(1)-\mathcal{L}^{n-1}(1)$, $\mathcal{L}^n(1)$ are pairwise
orthogonal and central in $\B_n$.  Hence the algebra  $\B_n$ is a
direct sum of ideals corresponding to these projections and
$i$-th  component of such a decomposition is isomorphic to
$(1-U^*U)\alpha^i(1)\A \alpha^i(1)$, if  $i=0,...,n-1$,  and
$\alpha^n(1)\A \alpha^n(1)$, if $i=n$. To see the latter it suffices to check that 
$$
U^iU^{*i}\A U^iU^{*i}= \alpha^i(1)\A\alpha^{i}(1)\ni a \to \mathcal{L}^i(a)=U^{*i}aU^i \in \B_n, \qquad i=1,...,n,
$$
is injective homomorphism which follows immediately from the fact that $U^i$ is a partial isometry. Accordingly, we get the isomorphism \eqref{coefficients izomor muchomor} and the proof is finished.

We note, cf. \cite[Prop. 2.2]{kwa2} or \cite[Prop 6.2]{kwa-leb1}, that 
$$
C^*(\A,U^*U)=U^*U\A \oplus (1-U^*U)\A\cong  \A/\ker\alpha \oplus \A/J 
$$
where $J$ is the ideal \eqref{ideal J}. This indicates that the  extended system $(\B,\alpha, \mathcal{L})$ can be reconstructed  from the triple $(\A,\alpha,J)$ and we will show that this is indeed the case. This will be achieved in section \ref{konstrukcja z uzuciem broni palnej}, but first we fix  notation and recall indispensable facts concerning transfer operators and covariant representations of  $C^*$-dynamical systems. 


\section{Transfer  operators and covariant representations}\label{transfer}
Let us fix   a $C^*$-dynamical system $(\A,\alpha)$.
A {\em transfer
operator} for $(\A,\alpha)$, see \cite{exel2},   is a positive linear map $\mathcal{L}:\A\to \A$ 
 such that
\begin{equation}
\mathcal{L}(\alpha(a)b) =a\mathcal{L}(b),\qquad a,b\in\A.
\label{b,,2}
\end{equation}
If additionally,  $\alpha (\mathcal{L}(1)) = \alpha (1)$ the transfer operator $\mathcal{L}$ is said to be \emph{non-degenerate} \cite{exel2}. The authors of  \cite{Bakht-Leb} called a transfer operator $\mathcal{L}$ for $(\A,\alpha)$ a \emph{complete transfer operator} if it satisfies 
\begin{equation}\label{definition of complete}
\alpha(\mathcal{L}(a))=\alpha(1)a \alpha(1), \qquad a\in \A.
\end{equation}
By  \cite{kwa3} it follows that a complete transfer operator exists if and only if  $\ker \alpha$ is   a complementary ideal in $\A$ and $\alpha(\A)$ is hereditary subalgebra of $\A$ (equivalently  $\ker\alpha$ is unital and $\alpha(\A)=\alpha(1)\A \alpha(1)$). 
 Then, see \cite{kwa3},  such a transfer operator is a unique non-degenerate transfer operator for $(\A,\alpha)$ and it is given by the formula
 $
 \mathcal{L}(a)=\alpha^{-1}(\alpha(1)a\alpha(1))$ 
where  $\alpha^{-1}$ is the inverse to the isomorphism $\alpha:(\ker\alpha)^\bot \to
\alpha(\A)$ and $(\ker\alpha)^\bot=\{a\in \A: a\ker\alpha=\{0\}\}$ is the \emph{annihilator} of $\ker\alpha$.
\begin{definition}[(cf. \cite{kwa-leb1})]\label{kowariant  rep defn*}
 A \emph{representation} of  $(\A,\alpha)$ is a  triple $(\pi,U,H)$
consisting of a unital faithful representation $\pi:\A\to L(H)$ in a
Hilbert space $H$ and an operator $U\in L(H)$ satisfying  
\begin{equation}\label{covariance rel1*}
U\pi(a)U^* =\pi(\alpha(a)),\qquad a \in \A.
\end{equation}
Then  $J=\{a\in \A: U^*U\pi(a)=\pi(a)\}$ is an ideal in $\A$ contained in $(\ker\alpha)^\bot$, cf. \cite[Cor. 1.5]{kwa-leb1} or \cite[Prop. 1.16]{kwa2}. We call $J$ the \emph{ideal of covariance} for $(\pi,U,H)$, and say that   $(\pi,U,H)$ is a \emph{$J$-covariant representation} of $(\A,\alpha)$. If $J=(\ker\alpha)^\bot$, we simply say that $(\pi,U,H)$ is a \emph{covariant representation} of $(\A,\alpha)$.
 \end{definition}
\begin{remark}
By \cite[Prop 1.10]{kwa-leb1} for each system $(\A,\alpha)$ and ideal $J$ in $(\ker\alpha)^\bot$ there exists a $J$-covariant representation of $(\A,\alpha)$.
\end{remark}

The next two statements explain to some extent the role of covariant representations (without prefix $J$) and  complete transfer operators.
\begin{proposition}\label{unital kowarians} 
Let $(\pi,U,H)$ be a representation of a $C^*$-dynamical system $(\A,\alpha)$ such that $\ker\alpha$ has a unit (is a complementary ideal in $\A$). The following conditions are equivalent:
\begin{itemize}
\item[i)] $(\pi,U,H)$ is a   covariant representation
\item[ii)] $U^*U\in \pi(\A)$
\item[iii)] $U^*U\in \pi(Z(\A))$ ($Z(\A)$ stands for the center of $\A$)
\item[iv)] $U^*U$ is the unit in $\pi((\ker\alpha)^\bot)$
\end{itemize}
In particular, if $\alpha$ is injective, then  $(\pi,U,H)$ is a   covariant representation if and only if $U$ is an isometry.
\end{proposition}

\textbf{Proof.} It is straightforward, as we know that $\pi(\ker\alpha)=(1-U^*U)\pi(\A)\cap \pi(\A)$ and $U^*U\in \pi(\A)'$, see \cite[Prop. 1.9]{kwa-leb1}.

\begin{proposition}\label{transfer covar rep}
Let $(\pi,U,H)$ be a representation of a $C^*$-dynamical system $(\A,\alpha)$  which admits 
a complete transfer operator  $\mathcal{L}:\A\to \A$. Then $(\pi,U,H)$ is a covariant representation if and only if 
\begin{equation}\label{transfer covariance}
\pi(\mathcal{L}(a))=U^*\pi(a)U, \qquad a\in \A.
\end{equation}
\end{proposition}
\textbf{Proof.} We recall, cf. \cite[Prop. 1.5]{kwa3}, that $\mathcal{L}(1)$ is the unit in $(\ker\alpha)^\bot$. Hence if $(\pi,U,H)$ is satisfies \eqref{transfer covariance},  then  $U^*U=\pi(\mathcal{L}(1))$ is the unit
in $\pi(\ker\alpha^\bot)$ and  $(\pi,U,H)$ is the covariant representation  of $(\A,\alpha)$ by   Proposition \ref{unital kowarians}.
Conversely, if $(\pi,U,H)$ is a covariant representation of $(\A,\alpha)$, then  $U^*U=\pi(\mathcal{L}(1)) \in \pi(Z(\A))$, again by  Proposition \ref{unital kowarians}, and  using \eqref{definition of complete} for $a\in \A$ we get
\begin{align*}
U^*\pi(a)U&=U^*(UU^*\pi(a) UU^*)U=U^*\pi( \alpha(1)a\alpha(1))U=U^*\pi(\alpha(\mathcal{L}(a)))U
\\
&=U^*U\pi(\mathcal{L}(a))U^*U=\pi(\mathcal{L}(1))\pi(\mathcal{L}(a))=\pi(\mathcal{L}(a)),
\end{align*}
which finishes the proof.

We  can always reduce investigation of $J$-covariant representations to covariant repre\-sentations (without prefix $J$) with the help of the following construction, cf. \cite[6.1]{kwa-leb1}, \cite[2.1.1]{kwa2}.
\begin{definition}\label{unitization definition}
Let $(\A,\alpha)$ be a $C^*$-dynamical system and let $J$ be an ideal in $(\ker\alpha)^\bot$. We treat $\A$ as a $C^*$-subalgebra of 
$$
\A_J=\big(\A/\ker\alpha\big) \oplus \big(\A/J\big)
$$
using the  embedding $
\A \ni a \longmapsto \big(a +\ker\alpha\big)\oplus \big(a + J\big)\in \A_J$.
We define an extension of  $\alpha$ up to $\A_J$, which we will still denote by $\alpha$,  by the formula
$$
\A_J\ni (a +\ker\alpha)\oplus (b + J)\stackrel{ }{\longrightarrow}(\alpha(a) +\ker\alpha)\oplus
(\alpha(a) + J)\in \A_J.
$$
We call   $(\A_J,\alpha)$ a  \emph{$C^*$-dynamical system obtained from $(\A,\alpha)$  by a $J$-unitization of the kernel}. 
\end{definition}
\begin{remark}
The kernel of  the endomorphism $\alpha:\A_J\to \A_J$ has the unit given by
$ (0 + \ker\alpha) \oplus (1 +J)$, and the algebras $\A_J$ and  $\A$ coincide if and only if 
$\ker\alpha$ is  unital and $J=(\ker\alpha)^\bot$. This  to some extent explains the terminology, cf. \cite[Rem. 6.1]{kwa-leb1}, \cite[Rem. 2.3]{kwa2}.
In the commutative case  passing from  $(\A,\alpha)$ to $(\A_J,\alpha)$ corresponds to compactification of the complement of the image of a partial mapping described in \cite[Prop. 2.4]{kwa2}.
\end{remark}

\begin{proposition}\label{unitization rep thm}
Let $(\A_J,\alpha)$ be a $C^*$-dynamical system obtained by a $J$-unitization of the kernel of $\alpha:\A\to \A$. There is a one-to-one correspondence between $J$-covariant
representations $ (\pi, U,H)$ of $(\A,\alpha)$ and covariant representations $(\pi_J, U,H)$ 
of $(\A_J,\alpha)$ established by the
equality
\begin{equation}\label{pi^+ formula}
\pi_J\big((a + \ker\alpha)\oplus (b + J)\big)=U^*U\pi(a) + (1-U^*U)\pi(b) .
\end{equation}
In particular, for every  $J$-covariant  representation  $(\pi, U,H)$ of $(\A,\alpha)$ the algebra
$
\A_J$ is isomorphic to $C^*\big(U^{*}U, \pi(\A)\big).
$
\end{proposition}
\textbf{Proof.} See \cite[Prop 6.2]{kwa-leb1} or \cite[Prop. 2.2]{kwa2}.


\section{Main construction}\label{konstrukcja z uzuciem broni palnej}


For convenience,   until Definition  \ref{naturalne rozszerzenie defn}
we assume that  the kernel of  $\alpha:\A \to \A$ is unital
and let $q$ denote the unit in $\ker \alpha$ (in general situation we will  pass through  the system $(\A_J,\alpha)$ described in Definition \ref{unitization definition}, see also Remark \ref{remark on q} below). We put
$$\label{algebry A_n 2}
\A_n:=\alpha^n(1)\A\alpha^n(1),\qquad n\in \bbN,
$$
 and define algebras $\B_n$ as direct sums of the form
 \begin{equation}\label{algebry B_n 2}
\B_n:=q\A_{0}\oplus q\A_{1}\oplus ... \oplus q\A_{n-1} \oplus \A_{n}, \qquad n\in \bbN.
\end{equation}
In particular, $\B_0= \A_0=\A$.  For each  $n\in \bbN$ we let $\alpha_n:\B_n \to \B_{n+1}$ to   be  a homomorphism schematically
presented by the diagram
$$
\begin{xy}
\xymatrix@C=3pt{
      **[r]  \B_n \ar[d]^{\alpha_n}& = &  q\A_{0} \ar[d]^{id} &  \oplus & ... & \oplus &
      q\A_{n-1}\ar[d]^{id}& \oplus & \A_{n} \ar[d]^{q}    \ar[rrd]^{\alpha}   \\
       \B_{n+1} & = &  q\A_{0}& \oplus & ... &  \oplus& q\A_{n-1} & \oplus & q\A_{n} & \oplus  & \A_{n+1}
        }
  \end{xy}
$$
and formally  given by  the formula
$$
\alpha_n (a_{0}\oplus ... \oplus a_{ n-1}\oplus a_{n})= a_{0}\oplus ... \oplus a_{ n-1}\oplus q a_{n}
\oplus  \alpha(a_{n}),
$$
where $a_{k}\in q\A_{k}$, $k=0,...,n-1,$ and $a_{n}\in \A_{n}$. Let us note that,   since $a_n=q a_n + (1-q) a_n$ and
$\alpha:(1-q)\A \to \alpha(\A)$ is an isomorphism, homomorphism  $\alpha_n$ is injective.
We define  $\B:=\underrightarrow{\lim\,\,}\{\B_n,\alpha_n\}$
to be the direct limit of the direct sequence \begin{equation}\label{ciag prosty C*-algebr}
\B_0\stackrel{\alpha_0}{\longrightarrow} \B_1
\stackrel{\alpha_1}{\longrightarrow}\B_2  \stackrel{\alpha_2}{\longrightarrow}...\,
\end{equation}
and  denote by $\phi_n: \B_n \to \B$, $n\in \bbN$, the natural
embeddings ($\phi_n$ are injective   since the bonding morphisms $\alpha_n$ are). Thus we have
$$
\phi_0(\A)=\phi_0(\B_0) \subset \phi_1(\B_1) \subset ...\subset \phi_n(\B_n)
\subset ...\qquad \textrm{ and } \qquad \B=\overline{\bigcup_{n\in\bbN} \phi_n(\B_n)}.
$$
We will identify the algebra $\A$ with the subalgebra $\phi_0(\A) \subset \B$ and
under this identification  we extend   $\alpha$ onto the algebra $\B$.
To this end, we  consider two sequences   (an inverse one and a direct one)
\begin{equation}\label{ciag odwrotny s}
\B_0\stackrel{s_1}{\longleftarrow} \B_1
\stackrel{s_2}{\longleftarrow}\B_2  \stackrel{s_3}{\longleftarrow}...\,,
\end{equation}
\begin{equation} \label{ciag prosty s_*}
\B_0
\stackrel{s_{*,0}}{\longrightarrow} \B_1
\stackrel{s_{*,1}}{\longrightarrow}\B_2  \stackrel{s_{*,2}}{\longrightarrow}...\,,
\end{equation}
where $s_{n}$ is a "left-shift" and $s_{*,n}$   is a  "right-shift":
$$
s_{n}(a_{0}\oplus a_{1}\oplus ... \oplus a_{{n}})=a_{1}\oplus a_{ 2}\oplus ... \oplus a_{n}
$$
$$
s_{*,n}(a_{0}\oplus ... \oplus a_{ n-1}\oplus a_{n})=0\oplus \big(\alpha(1)a_{0}\,
\alpha(1)\big)\oplus ...\oplus \big(\alpha^{n+1}(1)a_{n}\alpha^{n+1}(1)\big),
$$
 $a_{k}\in q\A_{k}$, $k=0,...,n-1$, $ a_{n}\in \A_{n}$.
Since $\alpha^{n}(1)$, $n\in \bbN$, form a  decreasing sequence of orthogonal projections,
 mappings $s_n$ and $s_{*,n}$ are well defined. Moreover the operators  $s_n$ are homomorphisms,
   whereas  operators $s_{*,n}$  in general fail to be multiplicative.
\begin{proposition}\label{granica prosta zepsute twierdzenie}
Sequence  \eqref{ciag odwrotny s} induces an endomorphism  $\alpha:\B\to \B$
extending the endomorphism $\alpha:\A\to \A$, whereas sequence \eqref{ciag prosty s_*}
induces an operator  $\mathcal{L}:\B\to \B$ which is a complete transfer operator for
 the extended $C^*$-dynamical system $(\B,\alpha)$.
\par
The word  'induces' means here that $\alpha$ and $\mathcal{L}$  are given on the dense
  $^*$-subalgebra $\bigcup_{n\in\bbN} \phi_n(\A_n)$ of $\B$ by  the formulae
\begin{equation}\label{wzory na delta i delta_*}
\alpha(a)= \phi_{n-1}(s_{n}(\phi_{n}^{-1}(a))), \qquad \mathcal{L}(a)= \phi_{n+1}(s_{*,n}(\phi_n^{-1}(a))),
\end{equation}
where  $a\in \phi_n(\B_n)$, $ n >0$.
\end{proposition}

\textbf{Proof.}
Direct computations show that the following  diagrams
$$
\begin{xy}
\xymatrix{
       \B_0 \ar[r]^{\alpha_0} & \B_1 \ar[r]^{\alpha_1} & \B_2 \ar[r]^{\alpha_2} &
       \,\,...\,\, \ar[r]^{\alpha_{n-1}}&  \B_{n} \ar[r]^{\alpha_n} & \,\, ...\,\,
       \\
  \B_0 \ar[r]^{\alpha_0} & \B_1 \ar[lu]_{s_1} \ar[r]^{\alpha_1} &\B_2 \ar[r]^{\alpha_2}
   \ar[lu]_{s_2}& \,\,...\,\, \ar[lu]_{s_{3}}\ar[r]^{\alpha_{n-1}}&  \B_{n} \ar[lu]_{s_{n}}\ar[r]^{\alpha_n} &
    \ar[lu]_{s_{n+1}}\,\, ...\,\,
              }
  \end{xy}
$$
$$
\begin{xy}
\xymatrix{
       \B_0 \ar[r]^{\alpha_0} & \B_1 \ar[r]^{\alpha_1} & \B_2 \ar[r]^{\alpha_2} & \,\,...\,\,
       \ar[r]^{\alpha_{n-1}}&  \B_{n} \ar[r]^{\alpha_n} & \,\, ...\,\,
       \\
  \B_0 \ar[ru]^{s_{*,0}} \ar[r]^{\alpha_0} & \B_1 \ar[ru]^{s_{*,1}} \ar[r]^{\alpha_1} &\B_2
  \ar[r]^{\alpha_2} \ar[ru]^{s_{*,2}}& \,\,...\,\,\ar[ru]^{s_{*,n-1}} \ar[r]^{\alpha_{n-1}}&
   \B_{n} \ar[ru]^{s_{*,n}}\ar[r]^{\alpha_n} & \,\, ...\,\,
              }
  \end{xy}
$$
commute. Hence \eqref{ciag odwrotny s} and \eqref{ciag prosty s_*} induce  certain linear
mappings on $\B$ (i.e.  formulae \eqref{wzory na delta i delta_*} make sense).
The former mapping, which for the sake of proof we denote by $\widetilde{\alpha}$, is a homomorphism (since    $s_n$ is a homomorphism for all $n\in \bbN$)
and the latter one, which we denote by $\mathcal{L}$, is positive   (because $s_{*n}$ posses that property for all $n\in \bbN$).
\\
We assert  that  the  mapping $\widetilde{\alpha}$ induced  by \eqref{ciag odwrotny s} agrees with $\alpha$ on  $\A$
which we identify with $\phi_0(\A)$. Indeed, an element $\phi_0(a)$, $a\in \A$,  of the inductive limit $\B$
is represented by the sequence $
(a, qa \oplus \alpha(a), qa \oplus q\alpha(a)\oplus  \alpha^2(a), ...)
$ and hence $\phi_{1}^{-1}(\phi_0(a))=qa \oplus \alpha(a)$. Thus in view of \eqref{wzory na delta i delta_*} we have
$$
\widetilde{\alpha}(\phi_0(a))= \phi_{0}(s_{1}(\phi_{1}^{-1}(\phi_0(a))))=\phi_{0}(s_{1}(qa \oplus \alpha(a)))=\phi_{0}( \alpha(a)).
$$
Thereby our assertion is true and we are justified to  denote by $\alpha$ the mapping $\widetilde{\alpha}$ induced by
\eqref{ciag odwrotny s}.
\\
To prove  that  $\mathcal{L}$ is   a complete transfer operator for $(\B,\alpha)$
 it suffices to show  \eqref{b,,2} and \eqref{definition of complete}.
 For that purpose we  take arbitrary elements $\a, \b \in \bigcup_{n\in\bbN} \phi_n(\B_n)\subset \B$
  and note that there exist   $n\in \bbN$,
  such that $\a=\phi_{n+1}(a)$ and $\b=\phi_{n}(b)$ for   $a\in  \B_{n+1}$ and  $b\in  \B_{n}$.
Direct computation shows that $s_{*,{n}}(s_{n+1}(a)b)=a\cdot s_{*,n}(b)$
  and thus using formulae \eqref{wzory na delta i delta_*} we have
\begin{align*}
\mathcal{L}(\alpha(\a)\b)&=\mathcal{L}(\phi_{n}(s_{n+1}(a))\,\b)=\phi_{n+1}(s_{*,n}(s_{n+1}(a)b))
\\
&=\phi_{n+1}(a\cdot  s_{*,n}(b))=\phi_{n+1}(a)\cdot \phi_{n+1}( s_{*,n}(b))=\a \mathcal{L}(\b)
\end{align*}
which proves \eqref{b,,2}. Similarly,  one checks
 that $s_{n+1}(s_{*,n}(a))=s_{n+1}(1)a
s_{n+1}(1)$ and then we have
\begin{align*}
\alpha(\mathcal{L}(\a))
&= \alpha(\phi_{n+1}(s_{*,n}(a))=\phi_{n}(s_{n+1}(s_{*,n}(a)))=\phi_{n}(s_{n+1}(1)a s_{n+1}(1) )
\\
&=\phi_{n}(s_{n+1}(1)) \phi_{n}(a) \phi_{n}(s_{n+1}(1))=\alpha(1) \a \alpha(1),
\end{align*}
which proves \eqref{definition of complete} and finishes the proof.

The systems $(\A,\alpha)$ and $(\B,\alpha)$ considered above coincide if and only if the range of $\alpha:\A\to \A$ is a hereditary subalgebra of  $\A$. Indeed,  the range of the endomorphism $\alpha:\B\to \B$ is always a hereditary subalgebra of $\B$,  as it admits a complete transfer operator. If $\alpha(\A)=\alpha(1)\alpha(1)$ is a hereditary subalgebra of  $\A$, then  $\A_n=\alpha^n(1)\A\alpha^n(1)=\alpha^n(\A)$, for all  $n\in \mathbb{N}$. Consequently,  the monomomorphisms $\alpha_n:\B_n \to \B_{n+1}$ are isomorphisms, and hence $(\A,\alpha)=(\B,\alpha)$, under our identifications.  
This justifies the following 

\begin{definition}\label{naturalne rozszerzenie defn}
If $(\A,\alpha)$ is such that $\ker\alpha$ is a complementary ideal, we call the  system $(\B,\alpha)$ described in Proposition \ref{granica prosta zepsute twierdzenie} a $C^*$-dynamical system obtained from $(\A,\alpha)$ by  \emph{hereditation of range}.
\end{definition}

\begin{theorem}\label{rozszerzenie a repr. kowariantne}
Suppose that  $(\B,\alpha)$ is a $C^*$-dynamical system obtained from $(\A,\alpha)$ by hereditation of the range.
There is a one-to-one correspondence between  covariant
representations $ (\pi, U,H)$ and $(\widetilde{\pi}, U,H)$ of
$(\A,\alpha)$ and $(\B,\alpha)$ respectively,  which is
established by the  relation
\begin{equation}\label{formula szpiegula}
 \widetilde{\pi}(\phi_n(a_{0}\oplus a_1 ... \oplus a_{n}))=
 \pi(a_{0})+ U^*\pi(a_{1})U +... + U^{*n}\pi(a_{n}) U^{n}.
\end{equation}
\end{theorem}

\textbf{Proof.} Let $(\widetilde{\pi}, U,H)$ be a  covariant representation   of $(\B,\alpha)$.
It is straightforward  that    $(\pi, U,H)$ where $\pi = \widetilde{\pi}|_{\A}$ is a
   representation of $(\A,\alpha)$. To see that  $(\pi, U,H)$ is a covariant representation, by Proposition \ref{unital kowarians}, it suffices to show that $q$ is the unit not only in the  kernel of $\alpha:\A\to \A$ but also in the kernel of its extension $\alpha:\B\to \B$.
To see the latter   let $n>0$ and notice that
$$
\phi_0(q)=\phi_n(q \oplus 0\oplus 0...\oplus 0).
$$
Thus for  $a=\phi_n(a_0\oplus a_1\oplus... \oplus a_n)\in \phi_n(\B_n)$  we have
$$
\alpha(a)=\phi_{n-1}(a_1\oplus... \oplus a_n)=0 \,\,\Longleftrightarrow\,\, a_1=...=a_n=0
\,\,\Longleftrightarrow\,\,a=qa.
$$
We fix now a  covariant representation  $(\pi, U,H)$   of $(\A,\alpha)$ and  show that  formula \eqref{formula szpiegula}
defines a faithful  representation $\widetilde{\pi}$ of $\B$. To this end,  we note that in view  of Proposition \ref{kowariantna podpowiedz do algebry wpsolczynnikow}
for every  $n\in\bbN$, the mapping $\widetilde{\pi}_n: \phi_n(\B_n) \to C^*\big(\bigcup_{k=0}^nU^{*n}\pi(\A) U^n \big)$
where
$$
 \widetilde{\pi}_n(\phi_n(a_{0}\oplus a_1 ... \oplus a_{n}))=
 \pi(a_{0})+ U^*\pi(a_{1})U +... + U^{*n}\pi(a_{n}) U^{n}
 $$
  is an isomorphism. Consequently, to show that   $\widetilde{\pi}: \B \to C^*\big(\bigcup_{n\in \bbN}U^{*n}\pi(\A) U^n \big)$ given by \eqref{formula szpiegula} is a well defined  isomorphism,  it suffices to check that the  diagram
  $$
  \begin{xy}
\xymatrix{ \B_n \ar[r]^{\alpha_n} \ar[d]_{\widetilde{\pi}_n \circ \phi_n} & \B_{n+1}\ar[d]^{\widetilde{\pi}_{n+1} \circ \phi_{n+1}}
       \\
   \widetilde{\pi}_n(\phi_n(\B_n)) \ar[r]^{id\quad\,\,} & \widetilde{\pi}_{n+1}(\phi_{n+1}(\B_{n+1}))
              }
  \end{xy}
  $$
commutes. Let $a=a_{0}\oplus a_1\oplus ... \oplus a_{n} \in \B_n$. Since $\pi(1-q)=U^*U$   we have
\begin{align*}
 U^{*n}\pi(a_n) U^{n} & = U^{*n}\pi(qa_n  + (1-q)a_n) U^{n}=  U^{*n}\pi(qa_n) U^{n}   + U^{*n}(U^*U)\pi(a_n)U^{n}
 \\
 & = U^{*n}\pi(qa_n) U^{n}   + U^{*n}(U^*U)\pi(a_n)(U^*U)U^{n}
 \\
 &= U^{*n}\pi(qa_n) U^{n}   + U^{*n+1}\pi(\alpha(a_n))U^{n+1}
\end{align*}
and thus
\begin{align*}
\widetilde{\pi}_n(\phi_n(a))&= \sum_{k=0}^{n-1}U^{*k}\pi(a_k) U^{k} + U^{*n}\pi(a_n) U^{n}  =  \widetilde{\pi}_{n+1}(\phi_{n+1}(\alpha_n (a))).
\end{align*}
 Accordingly,   $\widetilde{\pi}$ is a faithful representation of $\B$.  Since $U^*U\in \pi(\A)\subset\widetilde{\pi}(\B)$, in view of Proposition \ref{unital kowarians}, the only thing we need   to prove is that $(\widetilde{\pi}, U,H)$ is a representation of $(\B,\alpha)$. Let $a=\phi_n(a_0\oplus a_1\oplus... \oplus a_n)\in \phi_n(\B_n)$, $n
>0$.  Using the  relations $a_k\in \A_k=\alpha^k(1)\A \alpha^k(1)$, $a_0\in q\A=\ker\alpha$  and the  fact  that $\{U^{k*}U^{k}\}_{k=0}^\infty$ is a decreasing sequence of projections lying in the center of $ C^*\big(\bigcup_{n\in \bbN}U^{*n}\pi(\A) U^n \big)$, cf. \cite[Prop.  3.7]{Leb-Odz},  we have
\begin{align*}
U\widetilde{\pi}(a)U^*&= \sum_{k=0}^{n}U U^{*k}\pi(a_k) U^{k} U^*=\pi(\alpha(a_{0})) +  \sum_{k=1}^{n}UU^* U^{*k-1}\pi(a_k) U^{k-1}U U^*
\\
&=\sum_{k=1}^{n}UU^* (U^{*k-1}U^{k-1})U^{*k-1}\pi(a_k) U^{k-1} (U^{*k-1}U^{k-1})U U^*
\\
&=\sum_{k=1}^{n} (U^{*k-1}U^{k-1}) (UU^*) U^{*k-1}\pi(a_k) U^{k-1} (UU^*) (U^{*k-1}U^{k-1})
\\
&=\sum_{k=1}^{n} U^{*k-1} (U^{k} U^{*k})\pi(a_k) (U^{k} U^{*k})U^{k-1}
\\
& =\sum_{k=1}^{n} U^{*k-1} \pi(\alpha^k(1)a_k\alpha^{k}(1)) U^{k-1} =\sum_{k=0}^{n-1} U^{*k} \pi(a_{k+1}) U^{k}
\\
&=\widetilde{\pi}(\phi_{n-1}(a_1\oplus a_2\oplus...
\oplus a_n))=\widetilde{\pi}(\alpha(a)).
\end{align*}
This finishes the proof.

Putting together constructions from Definitions \ref{unitization definition} and \ref{naturalne rozszerzenie defn} we obtain a construction that embraces the general situation.
\begin{definition}\label{naturalne rozszerzenie defn}
Suppose  $(\A,\alpha)$ is an arbitrary  $C^*$-dynamical system and $J$ is an ideal in $\A$ such that $J\cap \ker\alpha=\{0\}$. Let $(\A_J,\alpha)$ be the  $C^*$-dynamical system obtained from $(\A,\alpha)$ by the $J$-unitization of the  kernel and let $(\B,\alpha)$ be the system obtained from $(\A_J,\alpha)$ by the hereditation of the range:
$$
\A\subset \A_J \subset \B.
$$
We  call the system $(\B,\alpha)$ a \emph{natural $J$-extension} of $(\A,\alpha)$.  If $J=(\ker\alpha)^\bot$ we  call $(\B,\alpha)$  simply a \emph{natural extension} of $(\A,\alpha)$.

\end{definition}
\begin{remark}\label{remark on q}
One can construct $(\B,\alpha)$ directly from $(\A,\alpha, J)$, without passing through $(\A_J,\alpha)$. To this end one may, almost literally, apply   our direct limit construction  changing only the meaning of $q$ from an element of $\A$ to the quotient map $q:\A\to \A/J$.
\end{remark}
\begin{remark} If $\ker\alpha$ is unital, the natural extension of $(\A,\alpha)$ coincides with the system obtained by the hereditation of the range of $\alpha$. 
The extended endomorphism 
$\alpha:\B \to \B$ is an automorphism if and only if $(\B,\alpha)$ is a natural extension of a unital monomorphism $\alpha:\A\to \A$ (to see it note that we always have $\alpha(\B)=\alpha(1)\B\alpha(1)$ and $1-\mathcal{L}(1)$ is the common unit in both of the kernels of $\alpha:\A_J\to\A_J$ and $\alpha:\B\to \B$).  
\end{remark}

In view of Proposition  \ref{unitization rep thm} and Theorem \ref{rozszerzenie a repr. kowariantne}, within the notation of the above definition, denoting  by $\mathcal{L}$ the complete transfer operator for $\alpha:\B\to \B$, we have 
$$
\A_J=C^*(\A,\mathcal{L}(1))=\mathcal{L}(1)\A\oplus (1-\mathcal{L}(1))\A,
$$
$$
\B=C^*\left(\bigcup_{n=0}^\infty \mathcal{L}^n(\A) \right)=\clsp\{\mathcal{L}^n(a): a\in \A, n\in \mathbb{N}\}.
$$
We also have the following statement.
\begin{theorem}\label{rozszerzenie a repr. kowariantne2}
Let $(\A,\alpha)$ be  an arbitrary $C^*$-dynamical system and let
$(\B,\alpha)$ be its natural $J$-extension with the complete transfer  operator
$\mathcal{L}$. There is a one-to-one correspondence between  $J$-covariant representations $ (\pi, U,H)$ of $(\A,\alpha)$ and
covariant representations $(\widetilde{\pi}, U,H)$ of  $(\B,\alpha)$,  which is established by the  relation
$$
 \widetilde{\pi}(\sum_{k=0}^n \mathcal{L}^k(a_k))=\sum_{k=0}^n U^{*k}\pi(a_{k})U^k,\quad a_k \in \A.
$$
In particular, 
\begin{itemize}
\item[i)] for every  $J$-covariant representation $(\pi, U,H)$ of $(\A,\alpha)$ we have
$$
\B \cong C^*\big(\bigcup_{n\in \bbN}U^{*n}\pi(\A) U^n \big)=\clsp\{U^{*n}\pi(a) U^n: a\in \A, n\in \bbN\};
$$
\item[ii)] the crossed product $C^*(\A,\alpha;J)$ of $\A$ by $\alpha$ associated to $J$ defined in \cite[Def. 1.12]{kwa-leb1} is naturally isomorphic to the crossed product $\B\rtimes_\alpha \bbZ$ defined in \cite[Def. 2.6]{Ant-Bakht-Leb} (cf. Proposition \ref{transfer covar rep}).
\end{itemize}
\end{theorem}

{\bf Acknowledgement.}  This work was in part supported by Polish  National Science Centre  grant number  DEC-2011/01/D/ST1/04112.



\end{document}